\documentclass[11pt,a4paper]{article}
\usepackage{amsmath}
\usepackage{amsfonts}
\usepackage{amssymb}
\usepackage{color}
\usepackage{epsfig}
\usepackage{graphics}
\newtheorem{lemma}{Lemma}

\newtheorem{conjecture}[lemma]{Conjecture}
\numberwithin{equation}{section} \numberwithin{lemma}{section}
\begin{document}

\title{\bf{Derivatives of random matrix characteristic polynomials with applications to elliptic curves }}
\vspace {2 in}

\author{N.\ C.\ Snaith\\
        School of Mathematics,\\ University of Bristol,\\
            Bristol BS8 1TW, UK}
\date{\today}
\maketitle \thispagestyle{empty} \vspace{.5cm}
\begin{abstract}
 The value distribution of derivatives of
 characteristic polynomials of matrices from $SO(N)$ is calculated
 at the point 1, the symmetry point on the unit circle of the
 eigenvalues of these matrices.  We consider subsets of matrices
 from $SO(N)$ that are constrained to have $n$ eigenvalues equal
 to 1, and investigate the first non-zero derivative of the characteristic polynomial at that point.
 The connection between the values of
random matrix characteristic polynomials and values of
$L$-functions in families has been well-established.  The
motivation for this work is the expectation that through this
connection with $L$-functions derived from families of elliptic
curves, and using the Birch and Swinnerton-Dyer conjecture to
relate values of the $L$-functions to the rank of elliptic curves,
 random matrix theory will be useful in probing important
questions concerning these ranks.
\end{abstract}

\section{Introduction}
\subsection{Random matrix theory and number theory}

The connection between random matrix theory and number theory
began with the work of Montgomery \cite{kn:mont73} when he
conjectured that the distribution of the complex zeros of the
Riemann zeta function follows the same statistics as the
eigenvalues of a random matrix chosen from $U(N)$ generated
uniformly with respect to Haar measure.  This conjecture is
supported by numerical evidence \cite{kn:odlyzko89} and also by
further work \cite{kn:hejhal94,kn:rudsar,kn:bogkea95,kn:bogkea96}
suggesting that the same conjecture is true for more general
$L$-functions. For all these $L$-functions there is a Generalized
Riemann Hypothesis that the non-trivial zeros lie on a vertical
line in the complex plane.  The conjectures mentioned above
concern the statistics of the zeros high on this critical line.

The philosophy of Katz and Sarnak
\cite{kn:katzsarnak99a,kn:katzsarnak99b} extended the connection
with random matrix theory by proposing that rather than averaging
over many zeros of a given $L$-function, if the zeros near to the
point where the critical line crosses the real axis are averaged
over a family of naturally connected $L$-functions, then they will
be found to follow the statistics of the eigenvalues of one of the
three classical compact groups of random matrices: $U(N)$, $O(N)$
or $USp(2N)$, where again the statistics are computed with respect
to the probability measure given by Haar measure. There is
numerical evidence for this conjecture as well \cite{kn:rub98},
and strong support is given to it by the rigorous work of Katz and
Sarnak in the case of function field zeta functions.

For a review of applications of random matrix theory to questions
in number theory see, for example, \cite{kn:conrey1} or
\cite{kn:keasna03}.

There has been a series of papers, starting with
\cite{kn:keasna00a} and continuing with
\cite{kn:confar00,kn:hughes00,kn:hughes03,kn:keasna00b,kn:cfkrs,kn:cfkrs1},
examining how random matrix theory can be used to predict the
distribution of values of the Riemann zeta function and other
$L$-functions, either averaged over an interval high on the
critical line, or over a family at the critical point where the
critical line crosses the real axis.  For large values of the
natural asymptotic parameter, for example the variable ordering
the $L$-functions within the family, the moments of $L$-functions
are conjectured to split into a product of an arithmetic
contribution, determined by the family being averaged over, and a
component derived from a random matrix calculation - the
corresponding moment of the characteristic polynomial of the
matrices in one of the three matrix groups $U(N)$, $O(N)$ or
$USp(2N)$. The asymptotic parameter on the random matrix side is
the dimension of the matrix $N$ and a natural equivalence can be
made between the two.

As will be reviewed in the following section, this conjecture for
the leading order behaviour of moments of $L$-functions has been
used by Conrey, Keating, Rubinstein and Snaith \cite{kn:ckrs00} to
predict the frequency, within a family of elliptic curves, of an
$L$-function taking the value zero at the critical point.  This is
done by using random matrix theory to predict the value
distribution of the $L$-values at the critical point and then the
discretization of these values
\cite{kn:waldspurger81,kn:shimura73,kn:kohzag81} to calculate a
probability of the $L$-value being zero.   See also David,
Fearnley and Kisilevsky \cite{kn:dfk04} for a similar use of
random matrix theory to predict frequency of vanishing at the
critical point amongst families of elliptic curve $L$-functions
twisted by cubic characters.

In Section \ref{sect:derivdist} we calculate the value
distribution of the first non-zero derivative of the
characteristic polynomial at the point one (corresponding to the
value at the critical point in the analogy with $L$-functions)
when a given number of eigenvalues are conditioned to lie at the
point one.  We find, at equation (\ref{eq:delaunay}), that moment
of the $n$th derivative grows like
\begin{equation} \label{eq:introdel}
 {\mathcal M}(n,M,s):=\langle |\Lambda_U^{(n)}(1)|^s
\rangle\sim(n!)^s (2\pi)^{s/2} 2^{-s^2/2-s(n-1)} \frac{G(n+1/2)}
{G(n+1/2+s)} M^{s^2/2+s(n-1/2)},
\end{equation}
where $\Lambda_U(e^{i\theta})$ is the characteristic polynomial of
$U$, $G(z)$ is the Barnes double gamma function (see \ref{eq:G})
and the angle brackets denote an average over the set of matrices
from $SO(N)$ with $n$ eigenvalues lying at the point one
($N=n+2M$). We also show that the probability that
$|\Lambda_U^{(n)}(1)|<X$ over those $U\in SO(N)$ with $n$
eigenvalues at 1 (again, $N=n+2M$) is, for small $X$, given by
(see (\ref{eq:nrzero}) and the sentence following)
\begin{equation}
 \tfrac{2}{2n+1} X^{\tfrac{2n+1}{2}}f(n,M),
\end{equation}
for the function $f$ given at (\ref{eq:f}).

If a discretization of the values of derivatives of $L$-functions
at the critical point were known, this would provide a way to
predict the frequency of zeros of various orders occurring at this
point. This has been investigated in the case of a single zero at
the critical point \cite{kn:crsw}.  In this situation, $n=1$, the
results of \cite{kn:crsw} give support to the validity of the
model presented here.

 The result presented in Section
\ref{sect:derivdist} on the moments of the first derivative of
characteristic polynomials from $SO(2N+1)$ (equation
\ref{eq:introdel} with $n=1$) has already been applied by Delaunay
in \cite{kn:delaunay05} in order to predict the moments of the
orders of Tate-Shafarevich groups and regulators of elliptic
curves with odd rank belonging to a family of quadratic twists.

In Section \ref{sect:evaluestats}, we briefly discuss the
eigenvalue statistics of random matrices with a certain number of
eigenvalues conditioned to lie at one.

We note that while we believe that the model presented in Section
\ref{sect:derivdist} should apply to the
 $L$-functions selected from families of quadratic twists of elliptic
curves for the property that they have high order zeros at the
critical point, there has been some very interesting theoretical
work computing the one- and two-level densities for parametric
families of elliptic curves that implies that the zeros of the
associated $L$-functions follow a different model (see
\cite{kn:mil04} and \cite{kn:young05}). This does not seem to be a
contradiction, as the zero statistics are examined over
collections of $L$-functions selected in very different ways, but
it certainly makes the question of the zero statistics of
$L$-functions with high-order zeros at the critical point very
intriguing.  See Section \ref{sect:discussion} for further
discussion on this issue.

\subsection{Random matrix theory and elliptic curves}

We review here the results of \cite{kn:ckrs00} which apply random
matrix theory to predicting the frequency of vanishing at the
critical point of the $L$-functions in the family of elliptic
curves described below; or equivalently, assuming the Birch and
Swinnerton-Dyer conjecture, the frequency of rank 2 or higher
curves occurring in the family of elliptic curves.  The motivation
for the random matrix calculations presented in Section
\ref{sect:derivdist} is the expectation that they may be used
similarly to examine curves of higher rank.

We consider an $L$-function (defined by the Dirichlet series and
Euler product below when ${\rm Re}\; s>3/2$)
\begin{equation}
L_{E}(s) = \sum_{n=1}^{\infty} \frac{a_n}{n^s}=\prod_{p} {\mathcal
L}_p(1/p^s)=\prod_{p\mid \Delta}
    \left(1-a_p p^{-s}\right)^{-1}
    \prod_{p\nmid \Delta}
    \left(1-a_p p^{-s}+p^{1-2s}\right)^{-1},
\end{equation}
that is associated to an elliptic curve $E$ over ${\mathbb{Q}}$
\begin{equation}
\label{eq:elliptic}
 E:\; y^2=x^3+Ax+B.
\end{equation}
The coefficients $a_p$, for prime $p$, are determined by
$a_p=p+1-\#E({\mathbb{F}}_p)$, where $\#E({\mathbb{F}}_p)$ counts
the number of pairs $x,y$, with $0\leq x,y\leq p-1$, such that
$y^2\equiv x^3+Ax+B\;(\mod \;p)$, plus one for the point at
infinity. $\Delta$ is the discriminant of the cubic $x^3+Ax+B$.
For an extremely clear introduction to elliptic curves in the
context discussed here, see the review paper by Rubin and
Silverberg \cite{kn:rubsil02}.

A family of elliptic curves is formed by
\begin{equation}
\label{eq:elliptic_d}
 E_d:\; dy^2=x^3+Ax+B
\end{equation}
for integers $d$ that are fundamental discriminants.  (A
fundamental discriminant is an integer other than 1 that is not
divisible by the square of any prime other than two and that
satisfies $d\equiv 1\mod 4$ or $d\equiv 8,12\mod 16$.)
 The corresponding family of $L$-functions, ordered by $|d|$, are
\begin{equation}
L_{E}(s,\chi_d)=\sum_{n=1}^{\infty} \frac{a_n\chi_d(n)}{n^s},
\end{equation}
where the characters $\chi_d(n)=\left(\frac{d}{n}\right)$ are
Kronecker's extension of the Legendre symbol, defined for prime
$p$ as
\begin{equation}
\left(\frac{d}{p}\right)=\left\{ \begin{array}{cc} +1& {\rm if}\;
    p\nmid d\; {\rm and}\;  x^2\equiv d({\rm mod} p)\; {\rm is\; soluble}\\
0& {\rm if}\; p|d \\ -1& {\rm if}\;
    p\nmid d\; {\rm and}\; x^2\equiv d({\rm mod} p) \; {\rm is\; not\; soluble}\end{array}\right..
\end{equation}
The family of curves $E_d$ is called the family of quadratic
twists of $E$ because the characters $\chi_d(n)$ are real,
quadratic Dirichlet characters.

Evidence (see for example \cite{kn:rub98}) points to the zeros
near the critical point (where the critical line crosses the real
axis) of such a family of $L$-functions having statistics like the
eigenvalues near 1 of matrices from $O(N)$ with Haar measure. To
explain this correspondence more clearly, about half of the
$L$-functions in this family are expected to have an even
functional equation that relates $L_E(s,\chi_d)$ to
$L_E(2-s,\chi_d)$, meaning that the zeros on the critical line
(${\rm Re} s=1$) appear in complex conjugate pairs and if there is
a zero at the critical point then it must be of even order, while
the other half have an odd functional equation, forcing a zero of
order 1,3,5,... at the critical point and all other zeros
appearing in complex conjugate pairs. Matrices in $O(N)$ also have
eigenvalues that appear in complex conjugate pairs, with the
possible exception of unpaired zeros at 1 and -1. The low-lying
zeros of elliptic curve $L$-functions from the family above that
have even functional equation show the statistics that eigenvalues
of $SO(2N)$ display near the point one.  These eigenvalues occur
in complex conjugate pairs and an eigenvalue at one must have even
multiplicity.  In contrast, the low-lying zeros of the
$L$-functions above that have an odd functional equation display
the same statistics as $SO(2N+1)$ eigenvalues near one, since in
this case there is always an eigenvalue at one and it has to have
odd multiplicity. (Note that we could just as well have chosen
$O^-(2N)$ matrices (orthogonal with determinant -1) to model the
$L$-functions with odd functional equation, since the statistics
of these eigenvalues near to one are the same as those of
$SO(2N+1)$, see for example \cite{kn:katzsarnak99a}.)

Both numerical and analytical evidence has been given already
\cite{kn:ckrs00} that random matrix theory can be used to
conjecture the frequency of $L$-functions with even functional
equation vanishing at the critical point in a family such as that
described above. This is particularly important because of the
Birch and Swinnerton-Dyer conjecture which suggests that the order
of the zero of an elliptic curve $L$-function at the critical
point is the same as the rank of the elliptic curve itself.  The
rank is a non-negative integer $r$ that characterizes the number
of rational points on the elliptic curve. The set of pairs $x,y\in
\mathbb{Q}$ that satisfy equation (\ref{eq:elliptic}), plus one
point at infinity, form a commutative group (with a standard law
of addition defined on the curve) that is isomorphic to
$\mathbb{Z}^r\bigoplus E(\mathbb{Q})_{tors}$.
$E(\mathbb{Q})_{tors}$ is the subgroup of elements of finite order
and $r$ is the rank of the curve.  The frequency of various ranks
in a family such as the one defined at (\ref{eq:elliptic_d}), or
even whether the rank is bounded in such a family of elliptic
curves, are important unanswered questions.

In \cite{kn:ckrs00} the first step is taken towards using random
matrix theory to generate conjectural answers to these questions
by addressing the frequency of $L$-functions from elliptic curve
families having any zero at all at the critical point.  The
argument used will be sketched below to illustrate the role played
by random matrix theory.  The purpose of this paper is to perform
the random matrix calculations needed to predict the frequency of
a zero of a {\it given order} at the critical point.  Applying the
random matrix results from Section \ref{sect:derivdist} to the
order of vanishing of $L$-functions involves the knowledge of
complicated arithmetic quantities and much more work is needed in
this area. Some preliminary numerical investigation has been done
\cite{kn:crsw} for the case of third order vanishing.

We now review the results in \cite{kn:ckrs00}.  We define the
family of $L$-functions discussed above with even functional
equation:
\begin{equation}
\label{eq:F+}
 \mathcal{F}_{E^+}=\{L_E(s,\chi_d) {\rm \; having\;
even\; functional \;equation}\}.
\end{equation}
The zeros of these $L$-functions near the critical point, $s=1$,
show the same statistics as the eigenvalues of matrices from
$SO(2N)$, and so the $L$-values at the critical point are modelled
by the characteristic polynomial
 \begin{equation}
 \label{eq:char}
\Lambda_U(e^{i\theta})=\prod_{n=1}^N \left(
1-e^{i(\theta_n-\theta)} \right) \left(1-
  e^{i(-\theta_n-\theta)} \right),
\end{equation}
evaluated at the point $\theta=0$.  Here $e^{\pm
i\theta_1},\ldots,e^{\pm i\theta_N}$ are the eigenvalues of the
matrix $U\in SO(2N)$.

The moments of $\Lambda_U(1)=\prod_{n=1}^N \big| 1-e^{i\theta_n}
\big|^2$ are easily calculated using Weyl's expression
\cite{kn:weyl} for Haar measure on the conjugacy classes of
$SO(2N)$ and a form of Selberg's integral to be
\cite{kn:keasna00b}

\begin{eqnarray}
&&\int_{SO(2N)}\Lambda_U(1)^sdU_{Haar}\nonumber \\
&&\;\;=2^{2Ns} \prod_{j=1}^N  \frac{\Gamma(N+j-1) \Gamma(s+j-1/2)}
{\Gamma(j-1/2) \Gamma(s+j+N-1)}\\
&&\equiv M_O(N,s).\nonumber
\end{eqnarray}

The $L$-function moments are then conjectured to have the form
\cite{kn:confar00,kn:keasna00b}

\begin{eqnarray}
\label{eq:elliptic_moment}
 M_E(T,s)&\equiv& \frac{1}{T^*}
\sum_{{|d|\leq T}\atop{L_E(s,\chi_d)\in \mathcal{F}_{E^+}}}
L_E(1,\chi_d)^s\nonumber\\
&\sim&a_s(E)M_O(N,s)
\end{eqnarray}
for large $T$.  Here $N=\log T$ (from equating the density of
zeros near the critical point with the density of the matrix
eigenvalues), the sum is over fundamental discriminants $d$, $T^*$
is the number of terms in the sum and $a_s(E)$ is an Euler product
that contains arithmetic information specific to the elliptic
curve $E$ and the family of $L$-functions being averaged over. In
practice, it is often a subset of $\mathcal{F}_{E^+}$ that is
summed over. If, for example, we select those $L$-functions
$L_E(s,\chi_d)$ in $\mathcal{F}_{E^+}$ with $d>0$ and further
restricted by a condition on $d \mod Q$, if $Q$ is odd, and on
$d\mod 4Q$, if $Q$ is even, then the arithmetic factor would be

\begin{eqnarray}
    a_s(E)= &&\prod_{p \nmid Q} \left( 1-p^{-1}\right)^{s(s-1)/2}
          \left(\frac{p}{p+1}\right)
          \left(
              \frac1p + \frac12
              \left(
                  \mathcal{L}_p(1/p)^{s} +
                  \mathcal{L}_p(-1/p)^{s}
              \right)
          \right)\nonumber \\
          \times&&\prod_{p \mid Q} \left( 1-p^{-1}\right)^{s(s-1)/2}
          \mathcal{L}_p(a_p/p)^{s}.
\end{eqnarray}
See \cite{kn:ckrs05} and \cite{kn:cfkrs}, Section 4.4, for more
examples.

Next we consider the distribution of the values of the
characteristic polynomials of $SO(2N)$ matrices at the point 1.
If $P_O(N,x)dx$ is the probability that the characteristic
polynomial of a matrix chosen from $SO(2N)$ with Haar measure has
a value between $x$ and $x+dx$, then
\begin{eqnarray}
\label{eq:dist}
P_O(N,x)&=&\frac{1}{2\pi ix} \int_{(c)} M_O(N,s)x^{-s}ds\nonumber\\
&\sim&x^{-1/2} h(N)
\end{eqnarray}
for $x\rightarrow 0^+$, since for small $x$ the behaviour is
dominated by the pole of $M_O(N,s)$ at $s=-1/2$.  Here $(c)$
denotes a path of integration along the vertical line from
$c-i\infty$ to $c+i\infty$.

For large $N$, $h(N)\sim 2^{-7/8}G(1/2)\pi^{-1/4} N^{3/8}$ ($G$ is
the Barnes double gamma function, defined as
\cite{kn:barnes00,kn:voros87}:
\begin{equation}
G(1+z)=(2\pi)^{z/2} e^{-[(1+\gamma)z^2+z]/2} \prod_{n=1}^{\infty}
\left[ (1+z/n)^n e^{-z+z^2/(2n)}\right],
\end{equation} where $\gamma$ is Euler's constant.  See also (\ref{eq:G}) for more properties of this function.)  Since the
probability that an element of $SO(2N)$ has a characteristic
polynomial whose value at 1 is $X$ or smaller is
$\int_0^XP_O(N,x)dx$, we find that the the behaviour of this
probability for small $X$ and large $N$ is
\begin{equation}
\lim_{N\rightarrow \infty}N^{-3/8} \lim_{X\rightarrow
0^+}\big(X^{-1/2} \int_0^XP_O(N,x)dx\big)=
2^{1/8}G(1/2)\pi^{-1/4}.
\end{equation}

We see from equation (\ref{eq:elliptic_moment}) that for large
$d$, moments of $L$-functions are conjectured to be just $a_s(E)$
(the prime product) times the random matrix moment $M_O(N,s)$. If
this is true, then we define $P_E(T,x)dx$ as the probability,
amongst members of $\mathcal{F}_{E^+}$, that $L_E(1,\chi_d)$, for
$|d|$ around $e^N$, will take a value between $x$ and $x+dx$,
giving

\begin{eqnarray}
 P_E(T,x) &=& \frac{1}{2\pi i x} \int_{(c)}
M_E(T,s)x^{-s} ds, \end{eqnarray} and an approximation for this
probability for small $x$ should be

\begin{eqnarray}
\label{eq:ellipticdist} P_E(T,x) &\sim& a_{-1/2}(E) x^{-1/2} h(N).
\end{eqnarray}
Here equating densities of zeros means $N\sim \log T$.

So, roughly, random matrix theory predicts that the probability
that an $L$-function $L_E(s,\chi_d)\in \mathcal{F}_{E^+}$, with
$|d|$ close to $t$, has a value at $s=1$ that is $X$ or smaller is
\begin{equation}
\label{eq:discret}
 \sim 2a_{-1/2}X^{1/2}2^{-7/8}G(1/2)\pi^{-1/4}
\log^{3/8}t.
\end{equation}

But these $L$-functions are constrained to take only certain
discretized values at the critical point $s=1$.  The $L$-values
have the form \cite{kn:waldspurger81,kn:shimura73,kn:kohzag81}:

\begin{equation}
L_E(1,\chi_d)=\kappa_E \frac{c_E(|d|)^2}{\sqrt{|d|}},
\end{equation}
where the $c_E(|d|)$ are integers, the Fourier coefficients of a
half-integral weight form.

So, if
\begin{equation}\label{eq:cutoff}L_E(1,\chi_d)<\frac{\kappa_E}{\sqrt{|d|}}\end{equation} then
\begin{equation}L_E(1,\chi_d)=0\end{equation} and the first thought would be to take
$X=\kappa_E /\sqrt{t}$ in (\ref{eq:discret}) and then to integrate
$t$ from 0 to $T$. However, there is arithmetical information
encoded in the $c_E(|d|)$ (one is not always the smallest allowed
non-zero value) and so to avoid this problem, the conjecture
stated in \cite{kn:ckrs00} is restricted to prime discriminants.

\begin{conjecture}(Conrey, Keating, Rubinstein, Snaith):
\label{conj:1}

Let $E$ be an elliptic curve defined over $\mathbb{Q} $.  Then
there is a constant $c_E\geq0$ such that

\begin{equation*}
\sum_{{p\leq T}\atop{{L_E(1,\chi_p)=0}\atop {L_E(s,\chi_p)\in
\mathcal{F}_{E^+}}}}\; 1\sim c_E T^{3/4} (\log T)^{-5/8}
\end{equation*}
\end{conjecture}
(The conjecture was originally stated in \cite{kn:ckrs00} with
$c_E>0$, but in fact numerics have revealed that the constant can
be zero \cite{kn:ckrs05} for certain families.  An explanation of
such a case with $c_E=0$ was given by Delaunay
\cite{kn:delaunay05a}.)

With the Birch and Swinnerton-Dyer conjecture, this suggests that
out of all the elliptic curves associated with $L$-functions in
$\mathcal{F}_{E^+}$ with prime discriminant $p\leq T$ (there are
of order $T/\log T$ of them), a number of order $T^{3/4}(\log
T)^{-5/8}$ of them should have rank two or greater.  The $T^{3/4}$
has been predicted previously by Sarnak using different arguments,
but random matrix theory adds more detailed information in the
form of the power on the logarithm. For numerical evidence
supporting the conjecture, see \cite{kn:ckrs00} and
\cite{kn:ckrs05}.

In the next section we will calculate the distribution of values
equivalent to (\ref{eq:dist}) that is required to probe questions
of elliptic curves of a given rank occurring in families of
quadratic twists (equation (\ref{eq:elliptic_d})).

\section{Random matrix calculations}

In this section we will calculate the probability of the $n$th
derivative of a characteristic polynomial (\ref{eq:cp}) of a
random $SO(N)$ matrix taking a value less then $X$ at $\theta=0$
when all lower derivatives are constrained to be zero at this
point.  That is, we perform the random matrix average over the
(measure zero) subset of $SO(N)$ of matrices with $n$ eigenvalues
at the point 1 on the unit circle, and $2M$ other eigenvalues
occurring in complex conjugate pairs ($N=n+2M$).

We will also illustrate the $n$-level density of the $2M$
eigenvalues {\em not} constrained to lie at 1.

\subsection{Value distribution of higher derivatives}
\label{sect:derivdist}

For a matrix $U\in SO(N)$ with $n$ eigenvalues equal to one
($N=n+2M$), the characteristic polynomial looks like
\begin{equation}
\label{eq:cp}
\Lambda_U(e^{i\theta})=(1-e^{-i\theta})^n\prod_{j=1}^M(1-e^{i(\theta_j-\theta)})
(1-e^{i(-\theta_j-\theta)}).
\end{equation}

We consider the $n$th derivative
\begin{eqnarray}
&&\frac{d^n}{d\alpha^n}(1-e^{-\alpha})^n\prod_{j=1}^M(1-e^{i\theta_j-\alpha})
(1-e^{-i\theta_j-\alpha})\big|_{\alpha=0}\nonumber\\
&&\qquad =\left[ \frac{d^n}{d\alpha^n}(1-e^{-\alpha})^n\right]
\prod_{j=1}^M
(1-e^{i\theta_j-\alpha})(1-e^{-i\theta_j-\alpha})\big|_{\alpha=0}\nonumber
\end{eqnarray}
and use this to define \begin{eqnarray}
 &&\Lambda_U^{(n)}(1)=n!
\prod_{j=1}^M (1-e^{i\theta_j})(1-
e^{-i\theta_j})=n!2^M\prod_{j=1}^M(1-\cos\theta_j).
\end{eqnarray}

The matrices with $n$ eigenvalues equal to one form a set of
measure zero within $SO(N)$ with Haar measure, for $n>1$.  For
even $N=2P$, Haar measure on the conjugacy classes of $SO(2P)$ is
given by
\begin{equation}
C' \prod_{1\leq j<k\leq P}(\cos\theta_k-\cos\theta_j)^2
d\theta_{1}\cdots d\theta_{P},
\end{equation}
with $C'$ a normalization constant, and we define the probability
of having $2q$ eigenvalues at 1 and the rest within infinitesimal
neighbourhoods of $e^{\pm i\theta_1},\ldots,e^{\pm i\theta_{P-q}}$
as
\begin{eqnarray}\label{eq:condprob}
&&\lim_{\alpha\rightarrow 0} \frac{\left(\int_{0}^{\alpha}\cdots
\int_0^{\alpha}\prod_{1\leq j<k\leq
P}(\cos\theta_k-\cos\theta_j)^2 d\theta_{P-q+1}\cdots
d\theta_{P}\right)d\theta_1\cdots
d\theta_{P-q}}{\int_0^{\pi}\cdots\int_0^{\pi}
\int_{0}^{\alpha}\cdots \int_0^{\alpha}\prod_{1\leq j<k\leq
P}(\cos\theta_k-\cos\theta_j)^2 d\theta_{P-q+1}\cdots
d\theta_{P}d\theta_1\cdots d\theta_{P-q}} \nonumber \\
&&=C(P,q)\prod_{j=1}^{P-q} (1-\cos\theta_j)^{2q} \prod_{1\leq
j<k\leq P-q} (\cos\theta_j-\cos\theta_k)^2d\theta_1\cdots
d\theta_{P-q}.
\end{eqnarray}
Here $C(P,q)$ is the normalization constant.

The same procedure for $N$ odd, starting from the measure on
$SO(2P+1)$
\begin{equation}
C''\prod_{j=1}^{P}(1-\cos\theta_j) \: \prod_{1\leq j<k\leq
P}(\cos\theta_k-\cos\theta_j)^2 d\theta_{1}\cdots d\theta_{P},
\end{equation}
with normalization constant $C''$, leads to a general expression
for the measure on the matrices in $SO(N)$ (where now $N$ can be
even or odd) with $n$ eigenvalues equal to 1 ($N=n+2M$):
\begin{equation}
\label{eq:measure} C(M,n)\prod_{j=1}^M (1-\cos\theta_j)^n
\prod_{1\leq j<k\leq M} (\cos\theta_j-\cos\theta_k)^2.
\end{equation}

Thus the quantity to calculate is
\begin{eqnarray}
\label{eq:integral} &&{\mathcal
M}(n,M,s):=\nonumber\\
&&C(M,n)\int_{0}^{\pi}\cdots \int_{0}^{\pi} |\Lambda_U^{(n)}(1)|^s
\prod_{j=1}^M(1-\cos\theta_j)^n\;\prod_{1\leq j<k\leq
M}(\cos\theta_j-\cos\theta_k)^2 d\theta_1\cdots d\theta_M\nonumber
\\
&&=C(M,n)\int_{0}^{\pi}\cdots \int_{0}^{\pi}(n!)^s 2^{Ms}
\prod_{j=1}^M (1-\cos\theta_j)^{n+s} \; \prod_{1\leq j<k\leq M}
(\cos\theta_j-\cos \theta_k)^2d\theta_1\cdots d\theta_M \nonumber
\\
&&=C(M,n) (n!)^s 2^{Ms} \int_{-1}^{1}\cdots \int_{-1}^{1}
\prod_{j=1}^M \frac{(1-x_j)^{n-1/2+s}}{(1+x_j)^{1/2}} \prod_{1\leq
j<k\leq M}(x_j-x_k)^2 dx_1\cdots dx_M.
\end{eqnarray}

This can be evaluated using a form of Selberg's integral (for
details see \cite{kn:mehta3}):

\begin{eqnarray}
\label{eq:Selberg} && \int_{-1}^1 \cdots \int_{-1} ^1 \prod_{1\leq
j< l\leq K} |(x_j-x_l)|^{2\gamma}  \prod_{j=1}^K
(1-x_j)^{\alpha-1}
(1+x_j)^{\beta-1} dx_j \nonumber \\
&&=2^{\gamma K(K-1)+K(\alpha+\beta-1)} \prod_{j=0}^{K-1}
\frac{\Gamma(1+\gamma+j\gamma)
  \Gamma(\alpha+j\gamma)\Gamma(\beta+j\gamma) }{ \Gamma(1+\gamma)
  \Gamma(\alpha+ \beta+\gamma(K+j-1))} ,
\end{eqnarray}

\noindent if ${\rm Re} \alpha>0$, ${\rm Re} \beta>0$ and ${\rm Re}
\gamma > -\min\left( \frac{1}{K}, \frac{{\rm Re} \alpha} {K-1},
\frac
  {{\rm Re} \beta} {K-1} \right)$.

We have $\gamma=1$, $\alpha=n+1/2+s$ and $\beta=1/2$, so the
integral in equation (\ref{eq:integral}) is
\begin{eqnarray}
&&{\mathcal M}(n,M,s)=C(M,n)(n!)^s 2^{2Ms+M(M-1)+Mn} \nonumber
\\
&&\qquad\qquad \times\prod_{j=1}^M
\frac{\Gamma(j+1)\Gamma(n-1/2+s+j)\Gamma(j-1/2)}
{\Gamma(n+M+s+j-1)}.
\end{eqnarray}
But the normalization constant can also be evaluated by Selberg's
integral to be $C(M,n)=\left(2^{M(M-1)+Mn}\prod_{j=1}^M
\frac{\Gamma(j+1)\Gamma(n-1/2+j)\Gamma(j-1/2)}
{\Gamma(n+M+j-1)}\right)^{-1}$, giving us
\begin{eqnarray}
{\mathcal M}(n,M,s)=(n!)^s 2^{2Ms} \prod_{j=1}^M
\frac{\Gamma(n-1/2+s+j)\Gamma(n+M+j-1)} {\Gamma(n-1/2+j)
\Gamma(n+s+M+j-1)}.
\end{eqnarray}

Consider $n$ fixed and finite and $M$ large.  If the probability
that $|\Lambda_U^{(n)}(1)|$ takes a value between $x$ and $x+dx$
is given by $P(n,M,x)dx$ then from a standard result in
probability (with $(c)$ denoting a path of integration along the
vertical line from $c-i\infty$ to $c+i\infty$)
\begin{eqnarray}
\label{eq:dists} P(n,M,x)&=& \frac{1}{2\pi ix} \int_{(c)}
x^{-s}{\mathcal M}(n,M,s) ds \\
&=& \frac{1}{2\pi ix} \int_{(c)} x^{-s} (n!)^{s}
2^{2Ms}\prod_{j=1}^M\frac{\Gamma(n-1/2+s+j)\Gamma(M+n+j-1)}
{\Gamma(n-1/2+j)\Gamma(n+s+M+j-1)} ds.\nonumber
\end{eqnarray}

\begin{figure}[htbp]
  \begin{center}
    \includegraphics[scale=0.72]
    {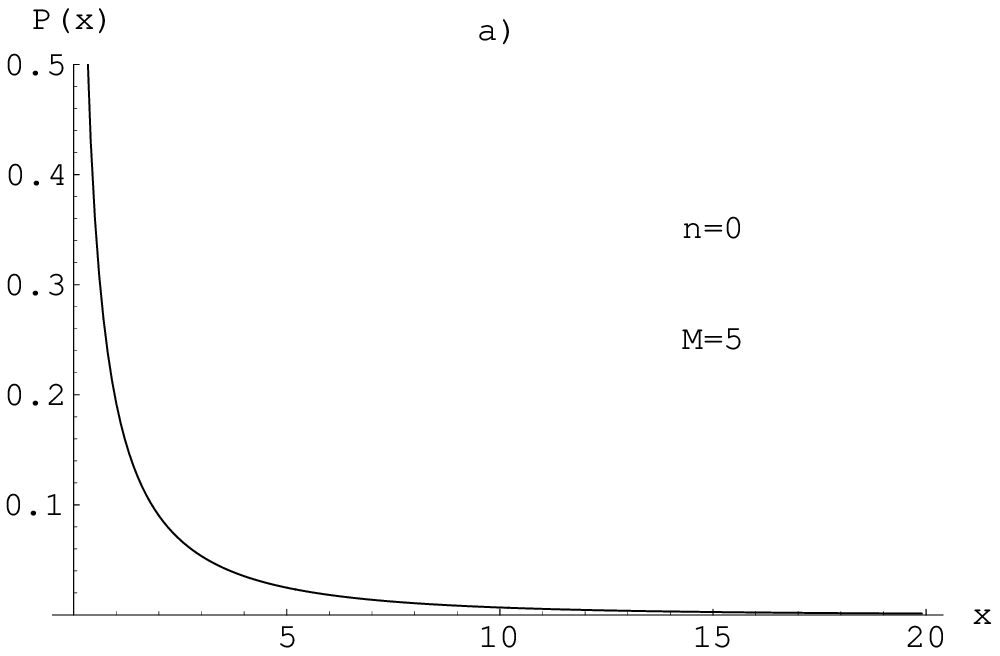}\\\vspace{0.1 in}
    \includegraphics[scale=0.72]{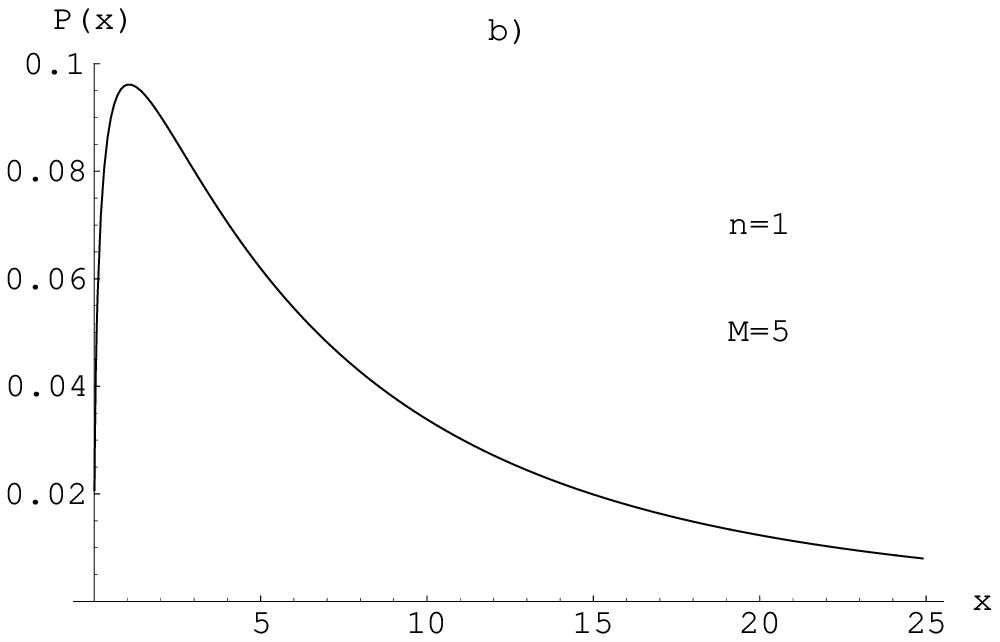}
\hspace{0.15 in}
\includegraphics[scale=0.72]
{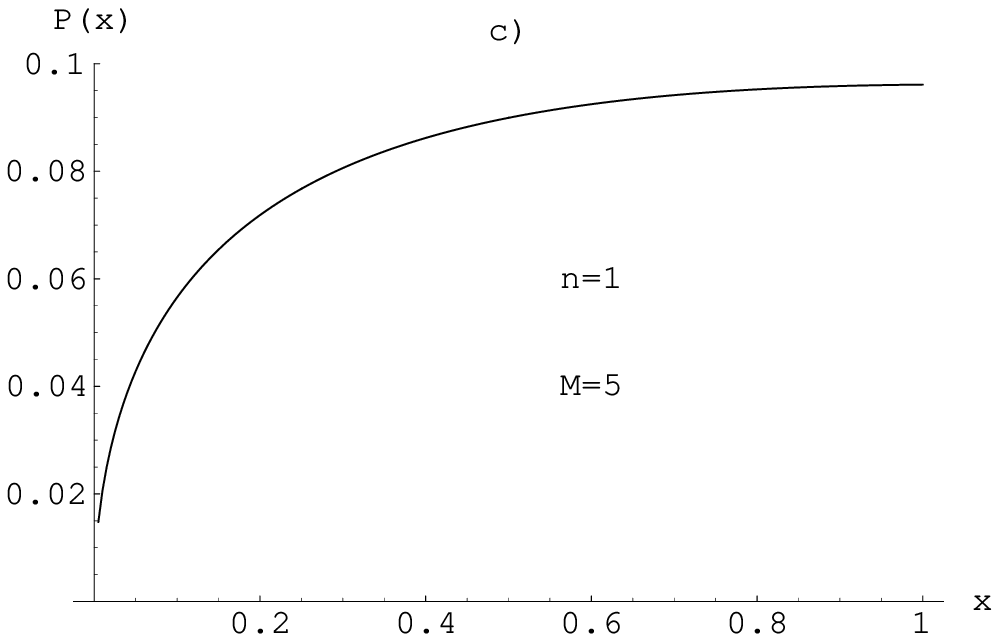}\\\vspace{0.1 in}
   \includegraphics[scale=0.72]{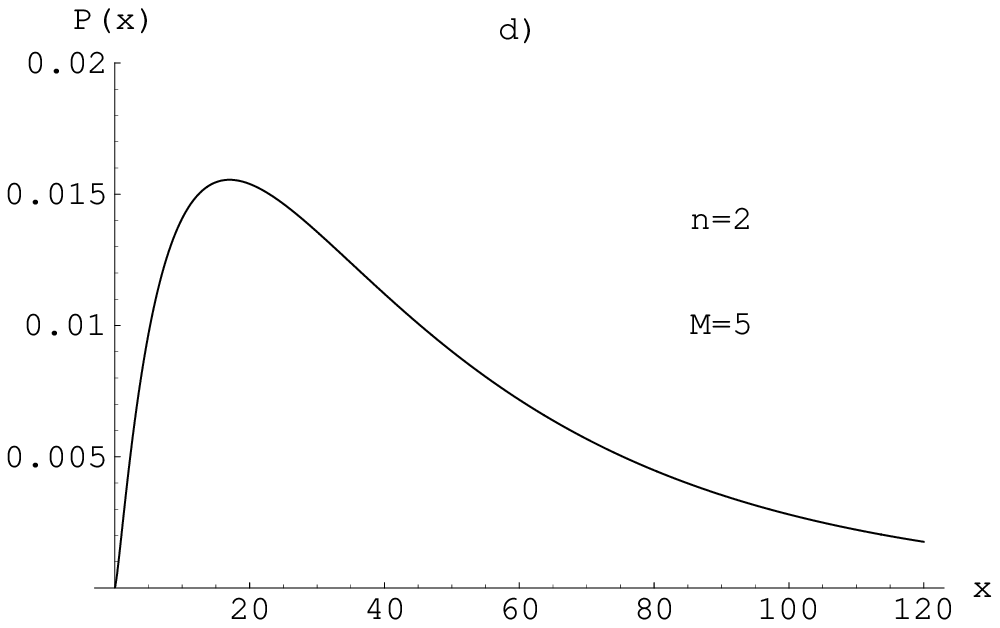}
\hspace{0.15 in}
\includegraphics[scale=0.72]
{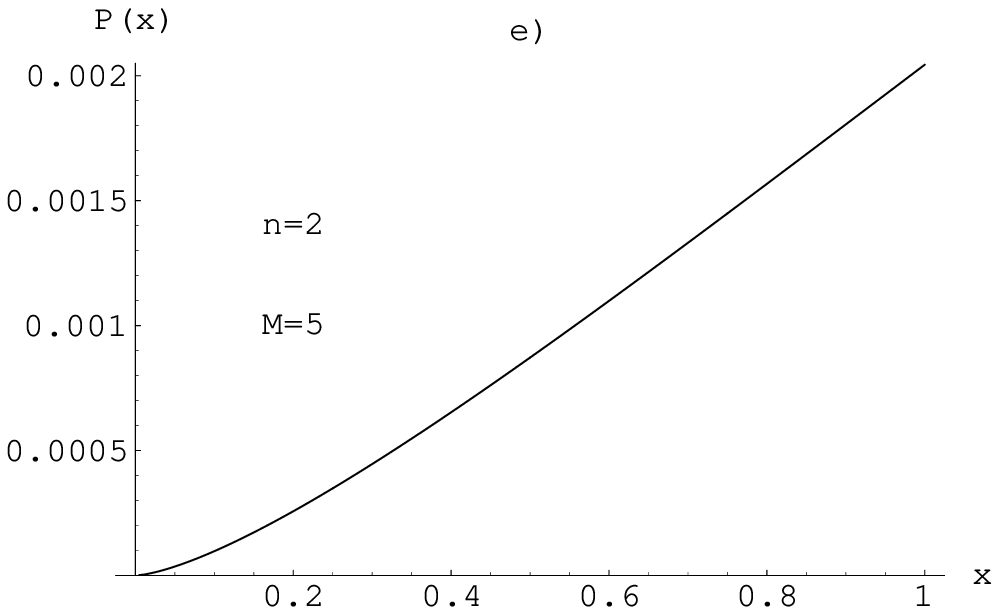}
    \caption{Distribution of values ($P(n,M,x)$ from (\ref{eq:dists}))
    at the point one of the $n$th derivative of the characteristic polynomial
    of the subset of matrices from $SO(n+2M)$ that are conditioned to
    have $n$ eigenvalues at one.  Figures c)
    and e) show in more detail the behaviour at the origin of figures b) and d)
    (see equation (\ref{eq:nrzero})).}
    \label{fig:dists}
  \end{center}
\end{figure}

We are interested in the behaviour at small $x$, and this is
dominated by the nearest pole to zero of the integrand: the pole
at $s=-(n+1/2)$ of $\Gamma(n-1/2+s+j)$. Thus for small $x$
\begin{equation}
\label{eq:nrzero} P(n,M,x)\sim x^{n-1/2}f(n,M).
\end{equation}
The probability that $|\Lambda_U^{(n)}(1)|<X$ over those $U\in
SO(N)$ with $n$ eigenvalues at 1 is therefore $\sim
\tfrac{2}{2n+1} X^{\tfrac{2n+1}{2}}f(n,M)$ for small $X$ and the
dependency on $n$ and $M$ is determined by
\begin{equation}\label{eq:f}
f(n,M)=(n!)^{\tfrac{-2n-1}{2}} 2^{-2M(\tfrac{2n+1}{2})}
\frac{1}{\Gamma(M)} \prod_{j=1}^M \frac{\Gamma(j)\Gamma(M+n+j-1)}
{\Gamma(n-1/2+j)\Gamma(M+j-3/2)}.
\end{equation}

We can discover the behaviour of $f(n,M)$ for large $M$ by
introducing the Barnes G-function \cite{kn:barnes00,kn:voros87}:
\begin{equation}
\label{eq:G} G(1+z)=(2\pi)^{z/2} e^{-[(1+\gamma)z^2+z]/2}
\prod_{n=1}^{\infty} \left[ (1+z/n)^n e^{-z+z^2/(2n)}\right],
\end{equation}

\noindent which has zeros at the negative integers, $-n$, with
multiplicity $n$ ($n=1,2,3 \ldots$).  Other properties useful to
us are that

\begin{eqnarray}
\label{eq:Grecurrence}
G(1)&=&1, \\
G(z+1)&=&\Gamma(z)\;G(z), \nonumber
\end{eqnarray}

\noindent and furthermore, for large $|z|$
\begin{eqnarray}
\label{eq:Gasymp} \log G(z+1)\sim z^2(\tfrac{1}{2}\log
z-\tfrac{3}{4})+\tfrac{1}{2}z\log (2\pi)-\tfrac{1}{12}\log z
+\zeta'(-1)+O(\tfrac{1}{z}).
\end{eqnarray}

Thus
\begin{equation}
\prod_{j=1}^{\ell} \Gamma(j+r)=\frac{G(\ell+r+1)}{G(1+r)}.
\end{equation}

So we can write
\begin{equation}
f(n,M)=(n!)^{\tfrac{-2n-1}{2}} 2^{-2M(\tfrac{2n+1}{2})}
\frac{G(M)G(2M+n)G(n+1/2)G(M-1/2)}{G(M+n)G(n+M+1/2)G(2M-1/2)}
\end{equation}
and expanding the $G$-functions for large $M$ gives
\begin{equation}
f(n,M)\sim (n!)^{\tfrac{-2n-1}{2}}
G(n+1/2)M^{-\tfrac{n^2}{2}+\tfrac{n}{2} +\tfrac{3}{8}}
2^{\tfrac{n^2}{2}-\tfrac{3n}{2}-\tfrac{7}{8}}
\pi^{-\tfrac{n}{2}-\tfrac{1}{4}}.
\end{equation}

Note that we can also use (\ref{eq:Gasymp}) to revisit ${\mathcal
M}(n,M,s)$ (defined at \ref{eq:integral}) and examine that
asymptotically for large $M$. This gives us the large $M$
behaviour of the moments at the point one of the $n$th derivative
of characteristic polynomials for which all lower derivatives
vanish.  We have
\begin{eqnarray}
\label{eq:delaunay}
{\mathcal M}(n,M,s)&=&(n!)^s 2^{2Ms} \frac{G(n+1/2+s+M)}
{G(n+1/2+s)} \frac{G(n+1/2)} {G(n+1/2+M)} \frac{G(n+2M)} {G(n+M)}
\frac{G(n+s+M)} {G(n+s+2M)}\nonumber \\
&\sim & (n!)^s (2\pi)^{s/2} 2^{-s^2/2-s(n-1)} \frac{G(n+1/2)}
{G(n+1/2+s)} M^{s^2/2+s(n-1/2)}.
\end{eqnarray}

\subsection{Eigenvalue statistics}
\label{sect:evaluestats}

To get an idea of what the eigenvalue statistics of the remaining
eigenvalues look like when we condition $n$ eigenvalues to lie at
one in the manner described in (\ref{eq:condprob}), we briefly
review the $m$-level densities of the unconditioned eigenvalues.
These calculations are not new; these same matrix ensembles were
studied by Nagao and Wadati \cite{kn:nagwad91} and also Duenez
\cite{kn:duenez04} and are discussed in the book by Forrester
\cite{kn:forrester}. Concurrently with the present work, this
matrix model has also been discussed in connection with the zero
statistics of elliptic curve $L$-functions by Miller
\cite{kn:mil05} (with appendix by Duenez). However, we include a
brief discussion of the statistics here to illustrate the effect
that forced eigenvalues at one have on the remaining eigenvalues.

The terminology for these statistics varies (for example, Mehta
\cite{kn:mehta3} refers to them as $m$-point correlation
functions) so we will write (with $C$ the normalization constant)
\begin{eqnarray}
&&C\int_{0}^{\pi}\cdots \int_{0}^{\pi} \sum_{1\leq j_1<\cdots<
j_m\leq N} f(\theta_{j_1},\ldots,\theta_{j_m})\;
(1-\cos\theta_j)^n\;\nonumber\\
&&\qquad \qquad\times\prod_{1\leq j<k\leq
M}(\cos\theta_j-\cos\theta_k)^2 d\theta_1\cdots d\theta_M\nonumber
\\
&&=\frac{1}{m!}\int_{0}^{\pi}\cdots
\int_{0}^{\pi}f(\theta_{1},\ldots,\theta_{m})\nonumber \\
&&\qquad
\qquad\times\det(K_M^{(n-1/2,-1/2)}(\theta_j,\theta_k))_{j,k=1,\ldots,m}\prod_{j=1}^m
(1-\cos\theta_j)^n d\theta_1\cdots d\theta_m\nonumber \\
&&=\frac{1}{m!}\int_{0}^{\pi}\cdots
\int_{0}^{\pi}f(\theta_{1},\ldots,\theta_{m})    \\
&&\qquad \qquad\times
\det(2^n\sin^n\tfrac{\theta_j}{2}\;\sin^n\tfrac{\theta_k} {2} \;
K_M^{(n-1/2,-1/2)}(\theta_j,\theta_k))_{j,k=1,\ldots,m}d\theta_1\cdots
d\theta_m    \nonumber
\end{eqnarray}
and refer to
$\det(2^n\sin^n\tfrac{\theta_j}{2}\;\sin^n\tfrac{\theta_k} {2} \;
K_M^{(n-1/2,-1/2)}(\theta_j,\theta_k))_{j,k=1,\ldots,m}$ as the
$m$-level density (see \cite{kn:conrey04} for further discussion
of these densities).    By standard techniques (see for example
\cite{kn:mehta3}, chapter 5, \cite{kn:szego}, chapter 4)
\begin{eqnarray}
&&K_M^{(n-1/2,-1/2)}(\theta_j,\theta_k))_{j,k=1,\ldots,m}\nonumber
\\
&&\quad=\sum_{j=0}^M(h_j^{(n-1/2,-1/2)}
)^{-1}P_j^{(n-1/2,-1/2)}(\cos\theta_j)\;P_j^{(n-1/2,-1/2)}(\cos\theta_k)\nonumber
\\
&&\quad= \frac{2^{-n+1}}{2M+n+1} \frac{\Gamma(M+2)\Gamma(M+n+1)}
{\Gamma(M+n+1/2)\Gamma(M+1/2)}\\
&& \times \nonumber\left(
\frac{P_{M+1}^{(n-1/2,-1/2)}(\cos\theta_j)\;P_M^{(n-1/2,-1/2)}(\cos\theta_k)-P_M^{(n-1/2,-1/2)}
(\cos\theta_j)\;P_{M+1}^{(n-1/2,-1/2)}(\cos\theta_k)}{\cos\theta_j-\cos\theta_k}\right)
\end{eqnarray}
is defined in terms of the Jacobi polynomials
$P_j^{(\alpha,\beta)}(\cos \theta)$ (described in detail in
\cite{kn:szego}) with the orthogonality condition
\begin{equation}
\int_{0}^{\pi} P_j^{(\alpha,\beta)}(\cos
\theta)\;P_k^{(\alpha,\beta)}(\cos\theta)(1-\cos
\theta)^nd\theta=\delta_{jk}h_m^{(\alpha,\beta)},
\end{equation}
where
\begin{equation}
h_j^{(\alpha,\beta)}=\frac{2^{\alpha+\beta+1}}{2j+\alpha+\beta+1}
\frac{\Gamma(j+\alpha+1)\Gamma(j+\beta+1)}
{\Gamma(j+1)\Gamma(j+\alpha+\beta+1)}.
\end{equation}

Scaling the eigenvalues and using asymptotic formulae for the
Jacobi polynomials (see \cite{kn:szego}, page 197) we find that in
the large-matrix limit (where $J_{\alpha}(x)$ are Bessel functions
of the first kind)
\begin{eqnarray}
L^{(n-1/2,-1/2)}(\theta,\phi)&:=&\lim_{M\rightarrow \infty}
\frac{\pi}{M} K_M^{(n-1/2,-1/2)}(\tfrac{\pi\theta}{M},\tfrac{\pi
\phi}{M}) 2^n\sin^n(\tfrac{\pi \theta}{2M})\sin^n(\tfrac{\pi
\phi}{2M}) \\
&=&(\pi\theta)^{1/2}(\pi\phi)^{1/2} \frac{\theta J_{n-3/2}(\theta
\pi)J_{n-1/2}(\phi \pi)-\phi J_{n-3/2}(\phi\pi)J_{n-1/2}(\theta
\pi)} {\phi^2-\theta^2}.\nonumber
\end{eqnarray}

\begin{figure}[htbp]
  \begin{center}
    \includegraphics[scale=1.5]
    {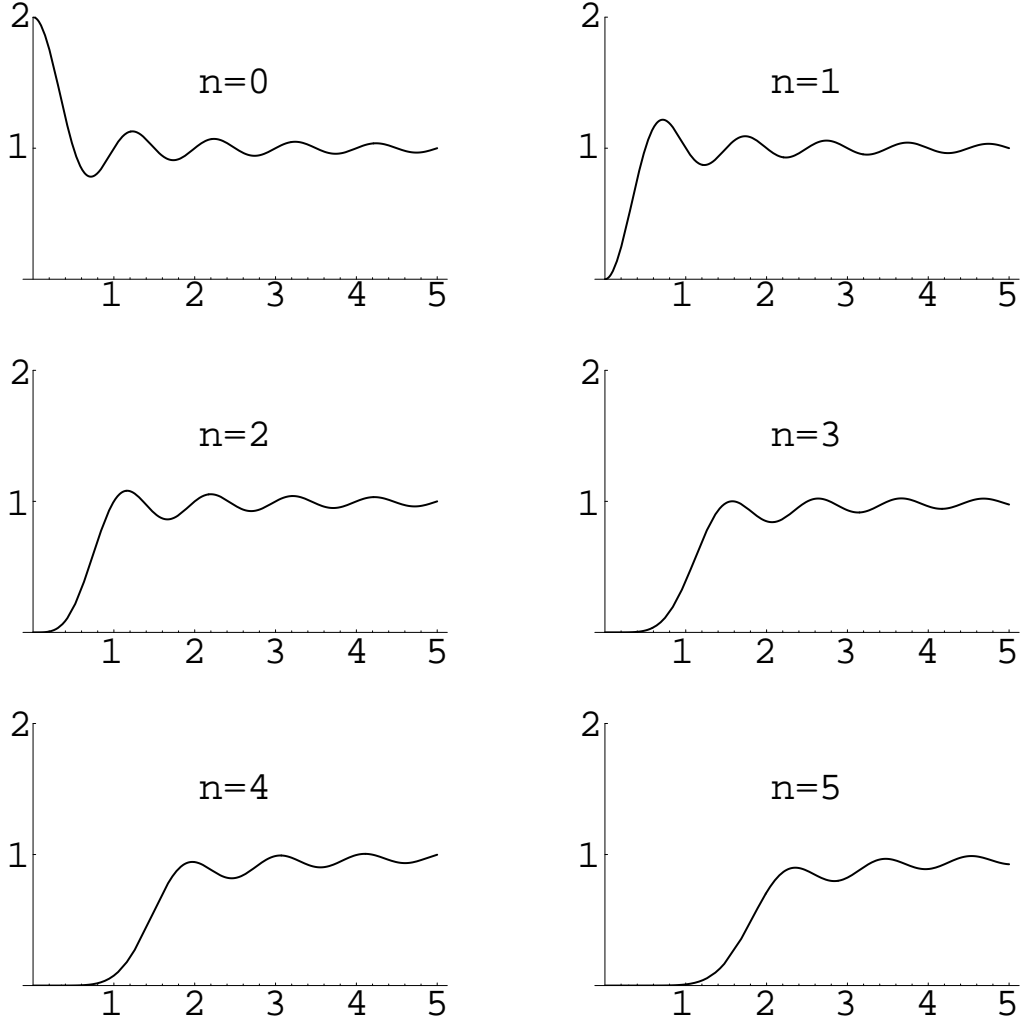}
    \caption{Scaled 1-level densities, from equation (\ref{eq:onelevel}) for various values of $n$, the degeneracy of the eigenvalue at one.}
    \label{fig:onedensity}
  \end{center}
\end{figure}

The scaled 1-level density is therefore
\begin{equation}
\label{eq:onelevel}
L^{(n-1/2,-1/2)}(\theta,\theta)=\frac{\pi^2}{2}
\theta(J_{n-3/2}^2(\theta\pi)+J_{n-1/2}^2(\theta \pi) -
\frac{2n-1}{\theta \pi} J_{n-1/2}(\theta
\pi)J_{n-3/2}(\theta\pi)),
\end{equation}
and can be seen to behave like a constant times $\theta^{2n}$ for
small $\theta$, where we remember that $n$ is the number of
degenerate eigenvalues. This increasing repulsion of the first
eigenvalue not located at 1 can be seen in Figure
\ref{fig:onedensity}.

\section{Discussion}
\label{sect:discussion}

We have shown that the probability that $|\Lambda_U^{(n)}(1)|<X$
over those $U\in SO(N)$ with $n$ eigenvalues at 1 ($N=n+2M$) is,
for small $X$, given by
\begin{equation}
\label{eq:disc1} \tfrac{2}{2n+1} X^{\tfrac{2n+1}{2}}f(n,M)
\end{equation}
and the mean value of the $n$th derivative grows like
\begin{equation}
\label{eq:del} {\mathcal M}(n,M,s):=\langle |\Lambda_U^{(n)}(1)|^s
\rangle\sim(n!)^s (2\pi)^{s/2} 2^{-s^2/2-s(n-1)} \frac{G(n+1/2)}
{G(n+1/2+s)} M^{s^2/2+s(n-1/2)},
\end{equation}
where the angle brackets denote an average over the set of
matrices from $SO(N)$ with $n$ eigenvalues lying at the point one
($N=n+2M$).

The result (\ref{eq:del}), see also (\ref{eq:delaunay}), was
applied by Delaunay in \cite{kn:delaunay05} to predict moments of
the orders of Tate-Shafarevich groups and the regulators of
elliptic curves belonging to a family of quadratic twists, see
(\ref{eq:elliptic_d}). The Birch and Swinnerton-Dyer conjecture
provides a formula for the first non-zero derivative of an
$L$-function in terms of various quantities related to the
associated elliptic curve.  One of these quantities is the order
of the Tate-Shafarevich group, and another, in the case of the
first derivative of $L$-functions with odd functional equation, is
the regulator.  For families of $L$-functions with even functional
equation the result (\ref{eq:delaunay}) is used with $n=0$ by
Delaunay to predict the asymptotic form of moments of the order of
the Tate-Shafarevich group for the associated family of elliptic
curves, and for families with odd functional equation the case
$n=1$ was used to conjecture the form of moments of the regulator.

To predict, for example, the number of $L$-functions in
$\mathcal{F}_{E^-}$ (the family defined like $\mathcal{F}_{E^+}$
in (\ref{eq:F+}) but with odd functional equation) that have a
zero of order at least three at the point 1 we need the result
(\ref{eq:disc1}) for $n=1$ plus some information analogous to
(\ref{eq:cutoff}) giving the smallest non-zero value that can be
taken by the derivative of an $L$-function from this family.  A
formula like (\ref{eq:cutoff}) is not known at the moment for
derivatives of $L$-functions, but numerics suggest \cite{kn:crsw}
that there is a gap between zero and the smallest non-zero value
of $L'$, which allows random matrix theory to be applied to the
question of the number of $L$-functions with a zero of order at
least three at the point 1.  Further work in this area is ongoing,
but it is also worth noting that the numerics in \cite{kn:crsw}
support the random matrix model presented here, as the numerical
cumulative distribution of values of $L'$ for $L$-functions in
$\mathcal{F}_{E^-}$ shows the expected $x^{3/2}$ behaviour near
$x=0$, agreeing with the exponent $3/2$ in (\ref{eq:disc1}) when
$n=1$.

This leads us to believe that the measure in (\ref{eq:measure})
should correctly model the zeros near the point 1 of $L$-functions
with at least an $n$th order zero at 1 selected from the family
$\mathcal{F}_{E^-}$ or $\mathcal{F}_{E^+}$ (depending on whether
$n$ is odd or even, respectively), but we note again that this is
not the only model for $L$-functions with $n$ zeros at 1.  The
measure presented here is what Miller \cite{kn:mil05} calls the
Interaction Model, and he contrasts this with the Independent
Model where the $n$ eigenvalues at 1 have no effect on the
statistics of the remaining eigenvalues.  This gives, for example
for matrices from $SO(2M+r)$ with $r$ eigenvalues at 1, the set of
matrices
$$\left\{\left(\begin{array}{cc}I_{r\times r}&
\\&g\end{array}\right): g\in SO(2M)\right\},$$ with a $r\times r$ identity block in the
upper left corner and with the joint probability density on the
remaining eigenvalues being Haar measure on $SO(2M)$, that is
$$\propto \prod_{1\leq j<k\leq M} (\cos \theta_k-\cos\theta_j)^2
\prod_{1\leq j\leq M}d\theta_j.$$  Miller finds numerical evidence
\cite{kn:mil05}, and for theoretical results see \cite{kn:mil04}
and \cite{kn:young05}, that the Independent Model agrees with the
zero statistics of $L$-functions of a parametric family of
elliptic curves constructed so that each member of the family has
rank $r$. Through the Birch and Swinnerton-Dyer conjecture, the
associated $L$-functions are expected to have a zero of order $r$
at the point 1, and Miller finds evidence for this high-order
zero, as well as for the fact that its presence does not affect
the statistics of the other zeros in the way that the model
presented here would (for example as illustrated in Figure 2), but
rather the zeros follow Miller's Independent Model. This is not a
contradiction to the proposal that the model presented here (the
same as Miller's Interaction Model) is correct for the family of
quadratic twists of elliptic curve $L$-functions.  It will be the
subject of further work to produce convincing numerics for the
zero statistics of the elliptic curve $L$-function families
discussed in this paper, although this is a difficult task due to
the small number of curves in these families with rank greater
than or equal to two.


\section{Acknowledgments}

Thanks to Brian Conrey for his number theoretical insight in
recognizing the potential importance of these calculations from
the start, to Peter Forrester for advice on orthogonal
polynomials, to Eduardo Due{\~ n}ez, Jon Keating and Steven Miller
for comments on early drafts of this paper and to Steve Gonek,
Chris Hughes, Michael Stoltz and Mark Watkins for asking and
answering a key question in the space of one lunchtime! My
gratitude also to the Royal Society for supporting my research
with the Dorothy Hodgkin Fellowship, to EPSRC for an Advanced
Research Fellowship and to the Isaac Newton Institute for
Mathematical Sciences for their hospitality during the programme
``Random Matrix Approaches in Number Theory" during which this
paper was prepared.



\end{document}